# Minimax Robust Function Reconstruction in Reproducing Kernel Hilbert Spaces

Richard J. Barton

***Abstract*— In this paper, the problem of approximating a function belonging to an arbitrary real-valued reproducing kernel Hilbert space (RKHS) on $\mathbb{R}^d$ based on *approximate* observations of the function is considered. The observations are approximate in the sense that the actual observations are known only to belong to a convex set of admissible observations. A minimax optimal approximation for the function is sought that minimizes the supremum of the RKHS norm on the error between the true function and the chosen approximation subject only to the conditions that the true function belongs to a uniformly bounded *uncertainty set* of functions that satisfy the constraints on the observations and that the approximation is a member of the RKHS. Such a solution is referred to as a *minimax robust reconstruction*. The solution to the minimax robust reconstruction problem is characterized and it is shown to be equivalent to solving a straightforward convex optimization problem. The stability properties of the minimax robust reconstruction are investigated and the approach is motivated by characterizing the minimax robust reconstruction for several specific convex observational models and discussing relationships with other approaches to function approximation.**

## I. Introduction

The results presented in this paper are related to the general problems of scattered data interpolation and approximation, and in particular to the large body of results on scattered data interpolation and approximation in a reproducing kernel Hilbert space (RKHS) setting [1-11]. We consider the problem of reconstructing an unknown function $f:\mathbb{R}^d \to \mathbb{R}$ based on *approximate* observations of the function on a subset $\mathcal{O}$ of points in $\mathbb{R}^d$. By approximate observations, we mean that the actual observations (that is, the true values of the function *f* on the set $\mathcal{O}$) are known only to belong to a convex set of admissible observations. In particular, we let[1] $f_{\mathcal{O}} \in \mathbb{R}^{\mathcal{O}}$ represent the observations (equivalently, the restriction of *f* to the set $\mathcal{O}$), and we assume that the observations are known only to satisfy $f_{\mathcal{O}} \in \mathcal{C}$, where $\mathcal{C}$ is an arbitrary convex subset of $\mathbb{R}^{\mathcal{O}}$. We further assume that $f \in \mathcal{H}$, where $\mathcal{H}$ is an RKHS with norm $\|\cdot\|_{\mathcal{H}}$, inner product $\langle \cdot,\cdot \rangle_{\mathcal{H}}$, and reproducing kernel $K:\mathbb{R}^d \otimes \mathbb{R}^d \to \mathbb{R}$, and that arbitrary functions in $\mathcal{H}$ are not determined uniquely[2] by their values on the set of points $\mathcal{O} \subset \mathbb{R}^d$. We seek a *minimax robust reconstruction* $\hat{f} \in \mathcal{H}$ of *f* that satisfies

$$\sup_{f\in\mathcal{U}} \left\| \hat{f} - f \right\|_{\mathcal{H}}^2 = \inf_{g\in\mathcal{H}} \sup_{f\in\mathcal{U}} \left\| g - f \right\|_{\mathcal{H}}^2, \tag{1}$$

where

$$\mathcal{U} = \left\{ f \in \mathcal{H} : \|f\|_{\mathcal{H}}^2 \le M,\ f_{\mathcal{O}} \in \mathcal{C} \right\} \ne \varnothing, \tag{2}$$

represents a non-empty *uncertainty set* of admissible reconstructions of $f$, and $M > 0$ is any real number satisfying

$$M \ge \inf_{f\in\mathcal{H}} \left\{ \|f\|_{\mathcal{H}}^2 : f_{\mathcal{O}} \in \mathcal{C} \right\}.$$

Note that putting some sort of bound on the norm for the set of admissible reconstructions is generally necessary in order to generate a well posed problem. That is, assuming that the set $\left\{ f \in \mathcal{H},\ f_{\mathcal{O}} \in \mathcal{C} \right\}$ is not empty and that a function $f \in \mathcal{H}$ is not uniquely determined by its values on the set $\mathcal{O}$, it is straightforward to show that for any function $g \in \mathcal{H}$, we have

$$\sup_{f\in\mathcal{H}} \left\{ \|g - f\|_{\mathcal{H}}^2 : f_{\mathcal{O}} \in \mathcal{C} \right\} = \infty .$$

Hence, without the bound on the norm, all approximations can be regarded as equally good (or equally bad) in terms of maximum possible deviation from the true function. On the other hand, we show below that the solution to Problem (1) does not actually depend on the value of *M* so that $\hat{f} \in \mathcal{H}$ solving (1) satisfies $\hat{f} = \lim_{M\to\infty} \hat{f}_M$, where

$$\sup_{f\in\mathcal{H}} \left\{ \left\| \hat{f}_M - f \right\|_{\mathcal{H}}^2 : \|f\|_{\mathcal{H}}^2 \le M,\ f_{\mathcal{O}} \in \mathcal{C} \right\}$$
$$= \inf_{g\in\mathcal{H}} \sup_{f\in\mathcal{H}} \left\{ \|g - f\|_{\mathcal{H}}^2 : \|f\|_{\mathcal{H}}^2 \le M,\ f_{\mathcal{O}} \in \mathcal{C} \right\}.$$

Since

[1] Here and elsewhere, the notation $\mathbb{R}^{\mathcal{O}}$ represents the set of functions $g:\mathcal{O} \to \mathbb{R}$.

[2] That is, we are only interested in the case when there are some degrees of freedom in the function reconstruction that are not determined by the set of observations.

$$\{ f \in \mathcal{H} : f_{\mathcal{O}} \in \mathcal{C} \} = \lim_{M\to\infty} \left\{ f \in \mathcal{H} : \| f \|_{\mathcal{H}}^2 \le M,\ f_{\mathcal{O}} \in \mathcal{C} \right\},$$

it follows that $\hat{f} \in \mathcal{H}$ solving (1) can be regarded as the unique meaningful solution to the problem

$$\sup_{f\in\mathcal{H}} \left\{ \| \hat{f} - f \|_{\mathcal{H}}^2 : f_{\mathcal{O}} \in \mathcal{C} \right\} = \inf_{g\in\mathcal{H}} \sup_{f\in\mathcal{H}} \left\{ \| g - f \|_{\mathcal{H}}^2 : f_{\mathcal{O}} \in \mathcal{C} \right\}.$$

Robust function reconstruction in contexts similar to the one considered here has been studied by many different investigators, and an excellent survey of those results is presented in [12, 13]. The particular approach taken here is of interest because of the very general form of the set $\mathcal{U}$ and because of the simple, intuitive form of the corresponding solution. Note that in the simplest case, the convex set of observations in (2) takes the form $\mathcal{C} = \{\mathbf{y}_0\}$, where $\mathbf{y}_0 \in \mathbb{R}^N$ is the nominal value of the function observed at $N$ points in $\mathbb{R}^d$. In this case, $\mathcal{C}$ consists of only a single point, and the uncertainty set $\mathcal{U}$ is just the set of all possible functions in $\mathcal{H}$ that interpolate between the observed values and satisfy an (essentially arbitrary) upper bound on the norm of *f*. As we discuss in the next section, the solution of (1) in this case is equivalent to the minimum-norm interpolation for *f*.

We show that with proper selection of the convex set $\mathcal{C}$, the solution of (1) can be identified with other well-known function approximation techniques such as smoothing splines and Moore-Penrose pseudo-inverses. Hence, part of the significance of this paper is that it presents a unified approach to function approximation that establishes for the first time the minimax optimality of function approximation techniques such as minimum-norm interpolation, function smoothing, and pseudo-inverses. These approaches have previously been justified and studied primarily on heuristic grounds. Furthermore, this approach provides considerable insight into stability and error analysis by stating and solving the approximation problem in a very geometrical context. Finally, from a practical point of view, the method easily accommodates a broad range of different observational models, accounts for the inevitable presence of noise in observations while still producing solutions that fall within a prescribed range of deviation from nominal behavior, and is computationally tractable.

The remainder of this paper is organized as follows. The basic results on minimax robust reconstruction are stated and proven in Section II. A discussion of computational complexity, stability characteristics and connections with other function approximation techniques is given in Section III. The method is illustrated by solving an approximation problem involving two-dimensional thin-plate splines in Section IV, and some concluding remarks are given in Section V.

## II. Results

Recall that in an RKHS $\mathcal{H}$ over $\mathbb{R}^d$ with kernel *K*, we have $K(\cdot, v) \in \mathcal{H}$ for any $v \in \mathbb{R}^d$. Furthermore, the reproducing property of the kernel *K* implies that

$$f(v) = \langle f, K(\cdot, v) \rangle_{\mathcal{H}},$$

for any $f \in \mathcal{H}$ and any $v \in \mathbb{R}^d$. Recall further that the set[3]

$$\mathcal{S} = \overline{\mathrm{sp}\{K(\cdot, u), u \in \mathcal{O}\}}, \tag{3}$$

consists of limits of functions of the form

$$s(v) = \sum_{i=1}^{n} a_i K(v, u_i),$$

for arbitrary $\{a_1, a_2, \ldots, a_n\} \subset \mathbb{R}$ and $\{u_1, u_2, \ldots, u_n\} \subset \mathcal{O}$, and that

$$\langle f, K(\cdot, u) \rangle_{\mathcal{H}} = \langle \mathcal{P}_{\mathcal{S}} f, K(\cdot, u) \rangle_{\mathcal{H}}, \quad \forall u \in \mathcal{O},$$

where $\mathcal{P}_{\mathcal{S}} f$ represents the projection of *f* onto $\mathcal{S}$. It follows from these relationships that for an uncertainty set $\mathcal{U}$ of the form (2), we have $f \in \mathcal{U}$ if and only if $\| f \|_{\mathcal{H}}^2 \le M$ and $\mathcal{P}_{\mathcal{S}} f \in \mathcal{U}$. Since we assume that a function $f \in \mathcal{H}$ is not uniquely determined by its values on the set $\mathcal{O}$, it also follows that $\mathcal{S} \ne \mathcal{H}$; that is $\mathcal{S}$ is a proper subset of $\mathcal{H}$ and so is $\bar{\mathcal{U}} \cap \mathcal{S}$. Our main result is as follows.

*Theorem 1.* Let $\mathcal{H}$ be an RKHS over $\mathbb{R}^d$ with reproducing kernel *K*, let $\mathcal{O} \subset \mathbb{R}^d$ be a (possibly infinite) set of observation points for an unknown function $f \in \mathcal{H}$, and assume that functions in $\mathcal{H}$ are not uniquely determined by their values on the set $\mathcal{O}$. Let $\mathcal{U}$ be given by (2), where $\mathcal{C} \subset \mathbb{R}^{\mathcal{O}}$ is convex, and let $\mathcal{S}$ be given by (3). Then a solution to Problem (1) always exists and is given by the unique element $\hat{f} \in \bar{\mathcal{U}} \cap \mathcal{S}$ satisfying

$$\| \hat{f} \|_{\mathcal{H}}^2 = \min_{f \in \bar{\mathcal{U}} \cap \mathcal{S}} \| f \|_{\mathcal{H}}^2 = \min_{f \in \bar{\mathcal{U}}} \| f \|_{\mathcal{H}}^2. \tag{4}$$

That is, the minimax robust reconstruction is the minimum norm element of $\bar{\mathcal{U}}$. It should be noted that the solution to (1) does not depend on the choice of the upper bound $M \ge \inf_{f\in\mathcal{H}} \left\{ \| f \|_{\mathcal{H}}^2 : f_{\mathcal{O}} \in \mathcal{C} \right\}$. To prove Theorem 1, we need the following two lemmas.

*Lemma 1.* $\hat{f} \in \mathcal{H}$ solves Problem (1) if and only if $\hat{f} \in \bar{\mathcal{U}} \cap \mathcal{S}$ and $\hat{f}$ satisfies

[3] Here and elsewhere, the over-bar notation $\bar{\mathcal{A}}$ indicates the closure of the set $\mathcal{A}$.

$$\sup_{f\in\mathcal{U}}\left\|\hat{f}-f\right\|_{\mathcal{H}}^2=\inf_{g\in\bar{\mathcal{U}}\cap\mathcal{S}}\sup_{f\in\mathcal{U}}\left\|g-f\right\|_{\mathcal{H}}^2;$$

that is, to find the optimal minimax reconstruction of the unknown function *f*, it is sufficient to consider functions $\hat{f}\in\bar{\mathcal{U}}\cap\mathcal{S}$.

*Proof.* Recall [14] that for any closed convex subset $\mathcal{K}\subset\mathcal{H}$ and any $g\in\mathcal{H}$, there is a unique element $g_0\in\mathcal{K}$, called the *projection of g onto* $\mathcal{K}$, such that

$$\|g-f\|_{\mathcal{H}}^2\ge\|g-g_0\|_{\mathcal{H}}^2,\quad\forall f\in\mathcal{K}.$$

Furthermore, $g_0\in\mathcal{K}$ satisfies

$$\langle f-g_0,g-g_0\rangle_{\mathcal{H}}\le 0,\quad\forall f\in\mathcal{K}.$$

Letting $\mathcal{K}=\bar{\mathcal{U}}$ and $g\in\mathcal{H}$, it follows that

$$\begin{aligned}\sup_{f\in\mathcal{U}}\|g-f\|_{\mathcal{H}}^2&=\sup_{f\in\mathcal{U}}\|g-g_0+g_0-f\|_{\mathcal{H}}^2\\&=\sup_{f\in\mathcal{U}}\left\{\begin{array}{l}\|g-g_0\|_{\mathcal{H}}^2-2\langle f-g_0,g-g_0\rangle_{\mathcal{H}}\\+\|g_0-f\|_{\mathcal{H}}^2\end{array}\right\}\\&\ge\sup_{f\in\mathcal{U}}\left\{\|g-g_0\|_{\mathcal{H}}^2+\|g_0-f\|_{\mathcal{H}}^2\right\}\\&\ge\sup_{f\in\mathcal{U}}\|g_0-f\|_{\mathcal{H}}^2\\&\ge\inf_{g\in\bar{\mathcal{U}}}\sup_{f\in\mathcal{U}}\|g-f\|_{\mathcal{H}}^2.\end{aligned}$$

Hence,

$$\inf_{g\in\mathcal{H}}\sup_{f\in\mathcal{U}}\|g-f\|_{\mathcal{H}}^2\ge\inf_{g\in\bar{\mathcal{U}}}\sup_{f\in\mathcal{U}}\|g-f\|_{\mathcal{H}}^2,$$

and since

$$\inf_{g\in\bar{\mathcal{U}}}\sup_{f\in\mathcal{U}}\|g-f\|_{\mathcal{H}}^2\ge\inf_{g\in\mathcal{H}}\sup_{f\in\mathcal{U}}\|g-f\|_{\mathcal{H}}^2,$$

we have

$$\inf_{g\in\bar{\mathcal{U}}}\sup_{f\in\mathcal{U}}\|g-f\|_{\mathcal{H}}^2\ge\inf_{g\in\mathcal{H}}\sup_{f\in\mathcal{U}}\|g-f\|_{\mathcal{H}}^2\ge\inf_{g\in\bar{\mathcal{U}}}\sup_{f\in\mathcal{U}}\|g-f\|_{\mathcal{H}}^2.$$

It follows that

$$\inf_{g\in\mathcal{H}}\sup_{f\in\mathcal{U}}\|g-f\|_{\mathcal{H}}^2=\inf_{g\in\bar{\mathcal{U}}}\sup_{f\in\mathcal{U}}\|g-f\|_{\mathcal{H}}^2.$$

Now, let $g\in\bar{\mathcal{U}}$ and $f\in\mathcal{U}$. Then $\mathcal{P}_{\mathcal{S}}g\in\bar{\mathcal{U}}\cap\mathcal{S}$ and

$$\begin{aligned}\|g-f\|_{\mathcal{H}}^2&=\|g-\mathcal{P}_{\mathcal{S}}g\|_{\mathcal{H}}^2+\|\mathcal{P}_{\mathcal{S}}g-\mathcal{P}_{\mathcal{S}}f\|_{\mathcal{H}}^2+\|\mathcal{P}_{\mathcal{S}}f-f\|_{\mathcal{H}}^2\\&\quad+2\langle g-\mathcal{P}_{\mathcal{S}}g,\mathcal{P}_{\mathcal{S}}f-f\rangle_{\mathcal{H}}\\&=\|g-\mathcal{P}_{\mathcal{S}}g\|_{\mathcal{H}}^2+\|\mathcal{P}_{\mathcal{S}}g-f\|_{\mathcal{H}}^2\\&\quad+2\langle g-\mathcal{P}_{\mathcal{S}}g,\mathcal{P}_{\mathcal{S}}f-f\rangle_{\mathcal{H}}.\end{aligned}$$

If $\langle g-\mathcal{P}_{\mathcal{S}}g,\mathcal{P}_{\mathcal{S}}f-f\rangle\ge 0$, then $\|g-f\|_{\mathcal{H}}^2\ge\|\mathcal{P}_{\mathcal{S}}g-f\|_{\mathcal{H}}^2$. If $\langle g-\mathcal{P}_{\mathcal{S}}g,\mathcal{P}_{\mathcal{S}}f-f\rangle<0$, let $\tilde{f}=\mathcal{P}_{\mathcal{S}}f-(f-\mathcal{P}_{\mathcal{S}}f)$. Then $\mathcal{P}_{\mathcal{S}}\tilde{f}=\mathcal{P}_{\mathcal{S}}f$, so $\tilde{f}\in\mathcal{U}$ and

$$\begin{aligned}\left\|g-\tilde{f}\right\|_{\mathcal{H}}^2&=\|g-\mathcal{P}_{\mathcal{S}}g\|_{\mathcal{H}}^2+\left\|\mathcal{P}_{\mathcal{S}}g-\mathcal{P}_{\mathcal{S}}\tilde{f}\right\|_{\mathcal{H}}^2+\left\|\mathcal{P}_{\mathcal{S}}\tilde{f}-f'\right\|_{\mathcal{H}}^2\\&\quad+2\left\langle g-\mathcal{P}_{\mathcal{S}}g,\mathcal{P}_{\mathcal{S}}\tilde{f}-\tilde{f}\right\rangle_{\mathcal{H}}\\&=\|g-\mathcal{P}_{\mathcal{S}}g\|_{\mathcal{H}}^2+\|\mathcal{P}_{\mathcal{S}}g-f\|_{\mathcal{H}}^2\\&\quad-2\langle g-\mathcal{P}_{\mathcal{S}}g,\mathcal{P}_{\mathcal{S}}f-f\rangle_{\mathcal{H}}\\&>\|\mathcal{P}_{\mathcal{S}}g-f\|_{\mathcal{H}}^2.\end{aligned}$$

It follows that,

$$\max\left\{\|g-f\|_{\mathcal{H}}^2,\left\|g-\tilde{f}\right\|_{\mathcal{H}}^2\right\}\ge\|\mathcal{P}_{\mathcal{S}}g-f\|_{\mathcal{H}}^2.$$

Hence, for all $f\in\mathcal{U}$, we have

$$\sup_{f\in\mathcal{U}}\|g-f\|_{\mathcal{H}}^2\ge\max\left\{\|g-f\|_{\mathcal{H}}^2,\left\|g-\tilde{f}\right\|_{\mathcal{H}}^2\right\}\ge\|\mathcal{P}_{\mathcal{S}}g-f\|_{\mathcal{H}}^2,$$

which implies that

$$\sup_{f\in\mathcal{U}}\|g-f\|_{\mathcal{H}}^2\ge\sup_{f\in\mathcal{U}}\|\mathcal{P}_{\mathcal{S}}g-f\|_{\mathcal{H}}^2\ge\inf_{g\in\bar{\mathcal{U}}\cap\mathcal{S}}\sup_{f\in\mathcal{U}}\|g-f\|_{\mathcal{H}}^2,$$

which in turn implies that

$$\begin{aligned}\inf_{g\in\bar{\mathcal{U}}\cap\mathcal{S}}\sup_{f\in\mathcal{U}}\|g-f\|_{\mathcal{H}}^2&\ge\inf_{g\in\bar{\mathcal{U}}}\sup_{f\in\mathcal{U}}\|g-f\|_{\mathcal{H}}^2\\&\ge\inf_{g\in\bar{\mathcal{U}}\cap\mathcal{S}}\sup_{f\in\mathcal{U}}\|g-f\|_{\mathcal{H}}^2.\end{aligned}$$

Hence,

$$\begin{aligned}\inf_{g\in\mathcal{H}}\sup_{f\in\mathcal{U}}\|g-f\|_{\mathcal{H}}^2&=\inf_{g\in\bar{\mathcal{U}}}\sup_{f\in\mathcal{U}}\|g-f\|_{\mathcal{H}}^2\\&=\inf_{g\in\bar{\mathcal{U}}\cap\mathcal{S}}\sup_{f\in\mathcal{U}}\|g-f\|_{\mathcal{H}}^2,\end{aligned}$$

which proves the lemma. ∎

*Lemma 2.* $\hat{f}\in\bar{\mathcal{U}}\cap\mathcal{S}$ solves Problem (1) if and only if

$$\begin{aligned}&\sup_{f\in\mathcal{U}}\left\{M-2\left\langle f,\hat{f}\right\rangle_{\mathcal{H}}+\left\|\hat{f}\right\|_{\mathcal{H}}^2\right\}\\&=\inf_{g\in\bar{\mathcal{U}}\cap\mathcal{S}}\sup_{f\in\mathcal{U}}\left\{M-2\langle f,g\rangle_{\mathcal{H}}+\|g\|_{\mathcal{H}}^2\right\}.\end{aligned}\tag{5}$$

*Proof.* Let $g\in\bar{\mathcal{U}}\cap\mathcal{S}$. Since, $\mathcal{S}\ne\mathcal{H}$, the orthogonal complement $\mathcal{S}^{\perp}$ of $\mathcal{S}$ in $\mathcal{H}$ contains non-zero elements, and for any $f\in\mathcal{U}$, $h\in\mathcal{S}^{\perp}\ne 0$, we have

$$\tilde{f}=\frac{\sqrt{M-\|\mathcal{P}_{\mathcal{S}}f\|_{\mathcal{H}}^2}}{\|h\|_{\mathcal{H}}}h+\mathcal{P}_{\mathcal{S}}f\in\mathcal{U},$$

with $\left\|\tilde{f}\right\|_{\mathcal{H}}^2=M$ and $\mathcal{P}_{\mathcal{S}}\tilde{f}=\mathcal{P}_{\mathcal{S}}f$. Hence,

$$
\begin{aligned}
\sup_{f\in\mathcal{U}} \|g-f\|_{\mathcal{H}}^2 &= \sup_{f\in\mathcal{U}} \left\{ \|f\|_{\mathcal{H}}^2 - 2\langle f,g\rangle_{\mathcal{H}} + \|g\|_{\mathcal{H}}^2 \right\} \\
&= \sup_{f\in\mathcal{U}} \left\{ \|f\|_{\mathcal{H}}^2 - 2\langle \mathcal{P}_{\mathcal{S}} f,g\rangle_{\mathcal{H}} + \|g\|_{\mathcal{H}}^2 \right\} \\
&= \sup_{f\in\mathcal{U}} \left\{ M - 2\langle \mathcal{P}_{\mathcal{S}} f,g\rangle_{\mathcal{H}} + \|g\|_{\mathcal{H}}^2 \right\} \\
&= \sup_{f\in\mathcal{U}} \left\{ M - 2\langle f,g\rangle_{\mathcal{H}} + \|g\|_{\mathcal{H}}^2 \right\}.
\end{aligned}
$$

It follows that,

$$
\begin{aligned}
&\inf_{g\in\bar{\mathcal{U}}\cap\mathcal{S}} \sup_{f\in\mathcal{U}} \|g-f\|_{\mathcal{H}}^2 \\
&= \inf_{g\in\bar{\mathcal{U}}\cap\mathcal{S}} \sup_{f\in\mathcal{U}} \left\{ M - 2\langle f,g\rangle_{\mathcal{H}} + \|g\|_{\mathcal{H}}^2 \right\}. \quad \blacksquare
\end{aligned}
$$

*Proof of Theorem 1*. The results of Lemma 1 and Lemma 2 imply that Theorem 1 will be proven if we can show that a solution to Problem (5) always exists and is given by the unique element $\hat{f} \in \bar{\mathcal{U}} \cap \mathcal{S}$ satisfying Equation (4). Toward that end, we define the *loss functional* $\mathcal{L} : \bar{\mathcal{U}} \times \left(\bar{\mathcal{U}} \cap \mathcal{S}\right) \to \mathbb{R}$ as

$$
\mathcal{L}(f,g) = M - 2\langle f,g\rangle_{\mathcal{H}} + \|g\|_{\mathcal{H}}^2 .
$$

To characterize the solutions to (5), we seek a *saddle point* for the *game* $\left(\bar{\mathcal{U}}, \bar{\mathcal{U}}\cap\mathcal{S}, \mathcal{L}\right)$; that is, we seek $(f_L, g_R) \in \bar{\mathcal{U}} \times \left(\bar{\mathcal{U}} \cap \mathcal{S}\right)$ such that

$$
\begin{aligned}
\mathcal{L}(f,g_R) &\le \mathcal{L}(f_L,g_R) \\
&\le \mathcal{L}(f_L,g), \quad \forall (f,g) \in \bar{\mathcal{U}} \times \left(\bar{\mathcal{U}} \cap \mathcal{S}\right).
\end{aligned} \tag{6}
$$

Since $\mathcal{L}$ is continuous, a solution to (6) satisfies

$$
\mathcal{L}(f_L,g_R) = \sup_{f\in\mathcal{U}} \mathcal{L}(f,g_R) \ge \inf_{g\in\bar{\mathcal{U}}\cap\mathcal{S}} \sup_{f\in\mathcal{U}} \mathcal{L}(f,g).
$$

Similarly, for any $g \in \bar{\mathcal{U}} \cap \mathcal{S}$,

$$
\mathcal{L}(f_L,g_R) \le \mathcal{L}(f_L,g) \le \sup_{f\in\mathcal{U}} \mathcal{L}(f,g).
$$

Hence,

$$
\mathcal{L}(f_L,g_R) \le \inf_{g\in\bar{\mathcal{U}}\cap\mathcal{S}} \sup_{f\in\mathcal{U}} \mathcal{L}(f,g),
$$

and it follows that

$$
\begin{aligned}
\mathcal{L}(f_L,g_R) &= \inf_{g\in\bar{\mathcal{U}}\cap\mathcal{S}} \sup_{f\in\mathcal{U}} \mathcal{L}(f,g) \\
&= \inf_{g\in\bar{\mathcal{U}}\cap\mathcal{S}} \sup_{f\in\mathcal{U}} \left( M - 2\langle f,g\rangle_{\mathcal{H}} + \|g\|_{\mathcal{H}}^2 \right).
\end{aligned}
$$

Substituting $\hat{f} = g_R$ gives the desired solution to (5). To find necessary and sufficient conditions for $(f_L, g_R) \in \bar{\mathcal{U}} \times \left(\bar{\mathcal{U}} \cap \mathcal{S}\right)$ to solve (6), we note that the Cauchy-Schwarz inequality implies

$$
\mathcal{L}(f_L,g) \ge M - \|\mathcal{P}_{\mathcal{S}} f_L\|_{\mathcal{H}}^2 = \mathcal{L}(f_L, \mathcal{P}_{\mathcal{S}} f_L)
$$

for all $g \in \bar{\mathcal{U}} \cap \mathcal{S}$, with equality if and only if $g = \mathcal{P}_{\mathcal{S}} f_L$. Hence,

$$
\mathcal{L}(f_L,g_R) \le \mathcal{L}(f_L,g), \quad \forall g \in \bar{\mathcal{U}} \cap \mathcal{S}
$$

if and only if $g_R = \mathcal{P}_{\mathcal{S}} f_L$. On the other hand, we have

$$
\mathcal{L}(f,g_R) \le \mathcal{L}(f_L,g_R), \quad \forall f \in \bar{\mathcal{U}}
$$

if and only if

$$
\begin{aligned}
&\lim_{\alpha\to 0} \frac{1}{\alpha}\left[ \mathcal{L}\left((1-\alpha) f_L + \alpha f, g_R\right) - \mathcal{L}(f_L,g_R) \right] \le 0 \\
&\Leftrightarrow \langle f,g_R\rangle_{\mathcal{H}} - \langle f_L,g_R\rangle_{\mathcal{H}} \ge 0 \\
&\Leftrightarrow \langle \mathcal{P}_{\mathcal{S}} f,g_R\rangle_{\mathcal{H}} - \langle f_L,g_R\rangle_{\mathcal{H}} \ge 0,
\end{aligned}
$$

where the notation "$\Leftrightarrow$" indicates that two statements are equivalent. Therefore, $(f_L, g_R) \in \bar{\mathcal{U}} \times \left(\bar{\mathcal{U}} \cap \mathcal{S}\right)$ solves (6) if and only if $g_R = \mathcal{P}_{\mathcal{S}} f_L$ and

$$
\begin{aligned}
&\langle \mathcal{P}_{\mathcal{S}} f,g_R\rangle_{\mathcal{H}} - \langle g_R,g_R\rangle_{\mathcal{H}} \ge 0, \quad \forall f \in \bar{\mathcal{U}} \\
&\Leftrightarrow \langle g,g_R\rangle_{\mathcal{H}} - \langle g_R,g_R\rangle_{\mathcal{H}} \ge 0, \quad \forall g \in \bar{\mathcal{U}} \cap \mathcal{S} \\
&\Leftrightarrow \langle g-g_R,g_R\rangle_{\mathcal{H}} \ge 0, \quad \forall g \in \bar{\mathcal{U}} \cap \mathcal{S} \\
&\Leftrightarrow \|g_R\|_{\mathcal{H}}^2 = \min_{g\in\bar{\mathcal{U}}\cap\mathcal{S}} \|g\|_{\mathcal{H}}^2,
\end{aligned}
$$

where the last statement follows from the definition of the projection of the origin onto $\bar{\mathcal{U}} \cap \mathcal{S}$ (*i.e.*, the minimum norm element of $\bar{\mathcal{U}} \cap \mathcal{S}$ ). This establishes the theorem. ■

## III. Discussion

It follows from Theorem 1 that solving Problem (1) is equivalent to finding the minimum-norm element of a convex set. From a practical point of view, one must also invert a possibly very large and ill-conditioned matrix, but the minimum-norm problem is still the most computationally complex part of the solution. Fortunately, finding the minimum-norm element of a convex set is a straightforward convex programming problem, and well known polynomial-time algorithms can be applied to find the solution [15-18]. To illustrate what is involved in solving a minimax robust reconstruction problem and to study the properties of the solutions, we characterize in this section the solutions for three particular convex uncertainty classes.

Throughout this section, we assume that the vector $\mathbf{y} = (y_0, y_1, \ldots, y_{N-1})^T$ represents a nominal vector of function observations at $N$ points $\mathcal{O} = \{v_0, v_1, \ldots, v_{N-1}\} \subset \mathbb{R}^d$, and we let the *gram matrix* for the problem be given by

$$
\mathbf{K} = \begin{bmatrix}
K(v_0,v_0) & K(v_0,v_1) & \cdots & K(v_0,v_{N-1}) \\
K(v_1,v_0) & K(v_1,v_1) & \cdots & K(v_1,v_{N-1}) \\
\vdots & \vdots & \cdots & \vdots \\
K(v_{N-1},v_0) & K(v_{N-1},v_1) & \cdots & K(v_{N-1},v_{N-1})
\end{bmatrix}.
$$

For simplicity, we assume that $\mathbf{K}$ is invertible. We focus on the three uncertainty sets

$$\mathcal{U}_i = \bar{\mathcal{U}}_i = \left\{ f \in \mathcal{H} \middle| \|f\|_{\mathcal{H}}^2 \le M, \mathbf{x} \in \mathcal{C}_i \right\}, \quad i = 1,2,3,$$

where

$$\mathbf{x} = (x_0, x_1, \ldots, x_{N-1})^T = \left( f(v_0), f(v_1), \ldots, f(v_{N-1}) \right)^T,$$

and

$$\mathcal{C}_1 = \left\{ \mathbf{x} \in \mathbb{R}^N \middle| \sum_{n=0}^{N-1} |x_n - y_n| \le \Delta \right\} = \mathcal{B}_1(\mathbf{y}, \Delta),$$

$$\mathcal{C}_2 = \left\{ \mathbf{x} \in \mathbb{R}^N \middle| \sum_{n=0}^{N-1} |x_n - y_n|^2 \le \Delta^2 \right\} = \mathcal{B}_2(\mathbf{y}, \Delta),$$

$$\mathcal{C}_3 = \left\{ \mathbf{x} \in \mathbb{R}^N \middle| \max_{n=0,1,\ldots,N-1} |x_n - y_n| \le \Delta \right\} = \mathcal{B}_\infty(\mathbf{y}, \Delta).$$

Here and elsewhere, the notation $\mathcal{B}_p(\mathbf{y}, \Delta)$ indicates a ball of radius $\Delta$ around the nominal observation vector $\mathbf{y}$ with respect to the $p$-norm in $\mathbb{R}^N$ for $p = 1, 2, \infty$. To avoid the trivial solution $\hat{f} = 0$, we assume that in each case, $\Delta > 0$ has been chosen such that $\mathcal{C}_i$ does not contain the origin. Then the subspace $\mathcal{S}$ consists of functions of the form

$$f(v) = \sum_{n=0}^{N-1} h_n K(v, v_n), \quad v \in \mathbb{R}^d,$$

where $\mathbf{h} = (h_0, h_1, \ldots, h_{N-1})^T$, and $f \in \bar{\mathcal{U}}_i \cap \mathcal{S}$ if and only if $f_{\mathcal{O}} = \mathbf{x} = \mathbf{K}\mathbf{h}$ for some $\mathbf{x} \in \mathcal{C}_i$, in which case we have

$$\|f\|_{\mathcal{H}}^2 = \mathbf{h}^T \mathbf{K} \mathbf{h} = \mathbf{x}^T \mathbf{K}^{-1} \mathbf{x} = \mathbf{h}^T \mathbf{x}.$$

Furthermore, if $f \in \bar{\mathcal{U}}_i \cap \mathcal{S}$ with $f_{\mathcal{O}} = \mathbf{x}_1 = \mathbf{K}\mathbf{h}_1 \in \mathcal{C}_i$ and $g \in \bar{\mathcal{U}}_i \cap \mathcal{S}$ with $g_{\mathcal{O}} = \mathbf{x}_2 = \mathbf{K}\mathbf{h}_2 \in \mathcal{C}_i$, then

$$\langle f, g \rangle_{\mathcal{H}} = \mathbf{h}_1^T \mathbf{K} \mathbf{h}_2 = \mathbf{x}_1^T \mathbf{K}^{-1} \mathbf{x}_2 = \mathbf{h}_1^T \mathbf{x}_2.$$

It follows from Theorem 1 that the minimax robust reconstruction $\hat{f} \in \bar{\mathcal{U}}_i \cap \mathcal{S}$ takes the form

$$\hat{f}(v) = \sum_{n=0}^{N-1} \hat{h}_n K(v, v_n), \quad v \in \mathbb{R}^d,$$

with $\hat{\mathbf{x}} = \mathbf{K}\hat{\mathbf{h}}$, where $\hat{\mathbf{x}} \in \mathcal{C}_i$ satisfies

$$\hat{\mathbf{h}}^T \hat{\mathbf{x}} = \left\| \hat{f} \right\|_{\mathcal{H}}^2 \le \|g\|_{\mathcal{H}}^2, \quad \forall g \in \bar{\mathcal{U}}_i \cap \mathcal{S}.$$

That is, $\hat{f}$ is the minimum-norm element of $\bar{\mathcal{U}}_i \cap \mathcal{S}$, or equivalently, $\hat{f}$ is the projection of the origin onto the closed convex set $\bar{\mathcal{U}}_i \cap \mathcal{S}$. Hence, for all $g \in \bar{\mathcal{U}}_i \cap \mathcal{S}$ of the form

$$g(v) = \sum_{n=0}^{N-1} h_n K(v, v_n), \quad v \in \mathbb{R}^d,$$

with $\mathbf{x} = \mathbf{K}\mathbf{h}$, we must have $\left\| \hat{f} \right\|_{\mathcal{H}}^2 \le \left\langle \hat{f}, g \right\rangle_{\mathcal{H}}$, or equivalently,

$$\hat{\mathbf{h}}^T \hat{\mathbf{x}} \le \hat{\mathbf{h}}^T \mathbf{x}, \quad \forall \mathbf{x} \in \mathcal{C}_i. \tag{7}$$

Hence, solving Problem (1) for each of the uncertainty sets $\mathcal{U}_i$, $i = 1,2,3$, is equivalent to identifying a pair $(\hat{\mathbf{h}}, \hat{\mathbf{x}})$ that satisfies inequality (7) with $\hat{\mathbf{x}} = \mathbf{K}\hat{\mathbf{h}}$. These three problems have been solved by Verdu and Poor [19] in the context of *robust matched filtering*. The solutions are discussed below.

$\mathcal{C}_1$: The pair $(\hat{\mathbf{h}}, \hat{\mathbf{x}})$ with $\hat{\mathbf{x}} = \mathbf{K}\hat{\mathbf{h}}$ satisfies Equation (7) if and only if

$$\hat{x}_n = y_n - \delta_n^2 \operatorname{sgn}(\hat{h}_n), \quad n = 0,1,\ldots,N-1,$$

where

$$\delta_n = 0 \text{ if } \left|\hat{h}_n\right| < \max_{n=0,1,,,\ldots N-1} \left|\hat{h}_n\right|,$$

$$\Delta = \sum_{n=1}^{N-1} \delta_n^2.$$

In the special case $\mathbf{K} = \operatorname{diag}(\sigma_0^2, \sigma_1^2, \ldots, \sigma_{N-1}^2)$, the *robust matched filter* $\hat{\mathbf{h}}$ is a clipped version of the *nominal matched filter* $\tilde{\mathbf{h}} = \mathbf{K}^{-1}\mathbf{y}$; that is,

$$\hat{h}_n = \begin{cases} \tilde{h}_n, & \text{if } \left|\tilde{h}_n\right| \le \tau, \\ \tau \operatorname{sgn}(\tilde{h}_n), & \text{if } \left|\tilde{h}_n\right| > \tau, \end{cases}$$

$$\Delta = \sum_{n=0}^{N-1} \sigma_n^2 \left( \left|\tilde{h}_n\right| - \tau \right)^+,$$

where $(x)^+ = \max(0, x)$.

$\mathcal{C}_2$: The pair $(\hat{\mathbf{h}}, \hat{\mathbf{x}})$ with $\hat{\mathbf{x}} = \mathbf{K}\hat{\mathbf{h}}$ satisfies Equation (7) if and only if

$$\hat{\mathbf{x}} = \mathbf{y} - \lambda \hat{\mathbf{h}},$$

$$\Delta = \lambda \left\| \hat{\mathbf{h}} \right\|_2.$$

where the symbol "$\|\cdot\|_2$" indicates the Euclidean norm of a vector. Alternatively,

$$\hat{\mathbf{h}} = (\mathbf{K} + \lambda \mathbf{I})^{-1} \mathbf{y},$$

$$\Delta = \lambda \left\| (\mathbf{K} + \lambda \mathbf{I})^{-1} \mathbf{y} \right\|_2.$$

$\mathcal{C}_3$: The pair $(\hat{\mathbf{h}}, \hat{\mathbf{x}})$ with $\hat{\mathbf{x}} = \mathbf{K}\hat{\mathbf{h}}$ satisfies Equation (7) if and only if

$$\hat{x}_n = \begin{cases} y_n - \Delta, & \text{if } \hat{h}_n > 0, \\ y_n + \Delta, & \text{if } \hat{h}_n < 0. \end{cases}$$

If $\tilde{\mathbf{h}} = \mathbf{K}^{-1}\mathbf{y}$ and

$$\Delta \le \frac{\min_i \left|\tilde{h}_i\right|}{\max_j \sum_{n=0}^{N-1} \left| \left(\mathbf{K}^{-1}\right)_{jn} \right|},$$

then

$$\hat{x}_n = \begin{cases} y_n - \Delta, & \text{if } \tilde{h}_n > 0, \\ y_n + \Delta, & \text{if } \tilde{h}_n < 0. \end{cases}$$

### A. *Stability and Connections with Other Approaches*

To get some insight into the stability of minimax robust reconstruction, we consider the solution for the uncertainty set $\mathcal{U}_2$ discussed above. In this case, the convex set of admissible observations is given by

$$\mathcal{C}_2 = \left\{ \mathbf{x} \in \mathbb{R}^N \,\middle|\, \sum_{n=0}^{N-1} \left| x_n - y_n \right|^2 \le \Delta^2 \right\} = \mathcal{B}_2\left(\mathbf{y}, \Delta\right).$$

Interestingly, the solution for Problem (1) in this case is identical to the solution of a classical RKHS smoothing problem [11, 20-25][4]. That is, the function $\hat{f}$ given by

$$\hat{f}\left(v\right) = \sum_{n=0}^{N-1} \hat{h}_n K\left(v, v_n\right), \quad v \in \mathbb{R}^d, \tag{8}$$

with

$$\hat{\mathbf{h}} = \left(\mathbf{K} + \lambda \mathbf{I}\right)^{-1} \mathbf{y}, \tag{9}$$

and

$$\Delta = \lambda \left\| \left(\mathbf{K} + \lambda \mathbf{I}\right)^{-1} \mathbf{y} \right\|_2 = \left\| \left(\tfrac{1}{\lambda}\mathbf{K} + \mathbf{I}\right)^{-1} \mathbf{y} \right\|_2, \tag{10}$$

also satisfies

$$\lambda \left\| \hat{f} \right\|_{\mathcal{H}}^2 + \sum_{n=1}^{N-1} \left| \hat{f}\left(v_n\right) - y_n \right|^2 = \inf_{f \in \mathcal{H}} \left\{ \lambda \left\| f \right\|_{\mathcal{H}}^2 + \sum_{n=1}^{N-1} \left| f\left(v_n\right) - y_n \right|^2 \right\}.$$

Hence, $\hat{f}$ can be regarded as a *generalized smoothing spline*, where the *smoothing parameter* $\lambda$ and the *uncertainty parameter* $\Delta$ are chosen jointly to satisfy (10). It follows that minimax robust reconstruction can be regarded as a generalization of function smoothing. Note that in the classical smoothing problem, one generally starts with a fixed value of $\lambda$, so that the uncertainty parameter becomes a function of the selected value of $\lambda$. In fact, Equation (10) implies that $\Delta \to 0$ as $\lambda \to 0$, $\Delta \to \left\| \mathbf{y} \right\|_2$ as $\lambda \to \infty$, and that, in general, the region of implied uncertainty increases with the norm of the nominal observed vector **y**. Note also that the equivalence established here between minimax robust reconstruction and function smoothing reveals the previously unrecognized minimax optimality of smoothing, which is generally justified and developed in a much more heuristic context.

As a result of the equivalence with smoothing in this particular case, one might conjecture that minimax robust reconstruction in general has the same stability properties as smoothing. That is, the general solution to the RKHS smoothing problem, given by (8) and (9) where $\lambda$ is a fixed parameter for all nominal observation vectors **y**, is unconditionally stable relative to variations in **y** since $\left(\mathbf{K} + \lambda \mathbf{I}\right)$ will not be ill-conditioned even if **K** is. Unfortunately, the stability properties of minimax robust reconstruction in general are not so well defined, and although the approach is much less sensitive to ill-conditioned gram matrices than techniques such as minimum-norm interpolation, it is not unconditionally stable even in the case illustrated here.

To see what is going on in some generality, let us assume that the set $\mathcal{C} \subset \mathbb{R}^N$ that determines the uncertainty set $\mathcal{U}$ is a closed convex set of small diameter with a fixed but arbitrary shape relative to its centroid **y**, which represents a nominal observation vector. In this case, a new nominal observation vector **y** changes the location but not the shape or the orientation of the set $\mathcal{C}$. While this is certainly not completely general, it covers a large class of useful uncertainty sets, and it is straightforward to visualize the behavior as **y** changes. The stability of minimax robust reconstruction in this case is determined by the rate of change in the solution $\hat{f}$ with respect to the norm $\left\| \cdot \right\|_{\mathcal{H}}$ relative to arbitrarily small changes in the coordinates of **y**. Furthermore, the quantity $\left\| \hat{f} \right\|_{\mathcal{H}}^2$ is given by

$$\left\| \hat{f} \right\|_{\mathcal{H}}^2 = \hat{\mathbf{x}}^T \mathbf{K}^{-1} \hat{\mathbf{x}},$$

where $\hat{\mathbf{x}}$ is the minimum-norm element of the set $\mathcal{C}$ with respect to the norm $\left\| \mathbf{x} \right\|_{\mathbf{K}}^2 = \mathbf{x}^T \mathbf{K}^{-1} \mathbf{x}$. To illustrate, we adopt the unit eigenvectors of **K** as our coordinate system, where $\left\{ \lambda_0 \ge \lambda_1 \ge \cdots \ge \lambda_{N-1} > 0 \right\}$ are the eigenvalues with associated orthonormal (with respect to $\left\| \cdot \right\|_2$) eigenvectors $\left\{ \mathbf{v}_0, \mathbf{v}_1, \ldots, \mathbf{v}_{N-1} \right\}$. The set $\mathcal{C}$ can be visualized (roughly) with respect to $\left\| \cdot \right\|_2$ as an *N*-dimensional sphere centered at **y**. On the other hand, with respect to $\left\| \cdot \right\|_{\mathbf{K}}$, which is really what matters, $\mathcal{C}$ can be visualized (again roughly) as an *N*-dimensional ellipsoid with major axis oriented along $\mathbf{v}_{N-1}$ and minor axis oriented along $\mathbf{v}_0$. Clearly, for visualization purposes, we have reverted to the earlier example with $\mathcal{C} = \mathcal{C}_2$. This is illustrated in Figure 1 for the two-dimensional subspace spanned by $\mathbf{v}_0$ and $\mathbf{v}_{N-1}$.

[4] This is a generalization of a similar well-known relationship for univariate function approximation with polynomial splines [24, 25].

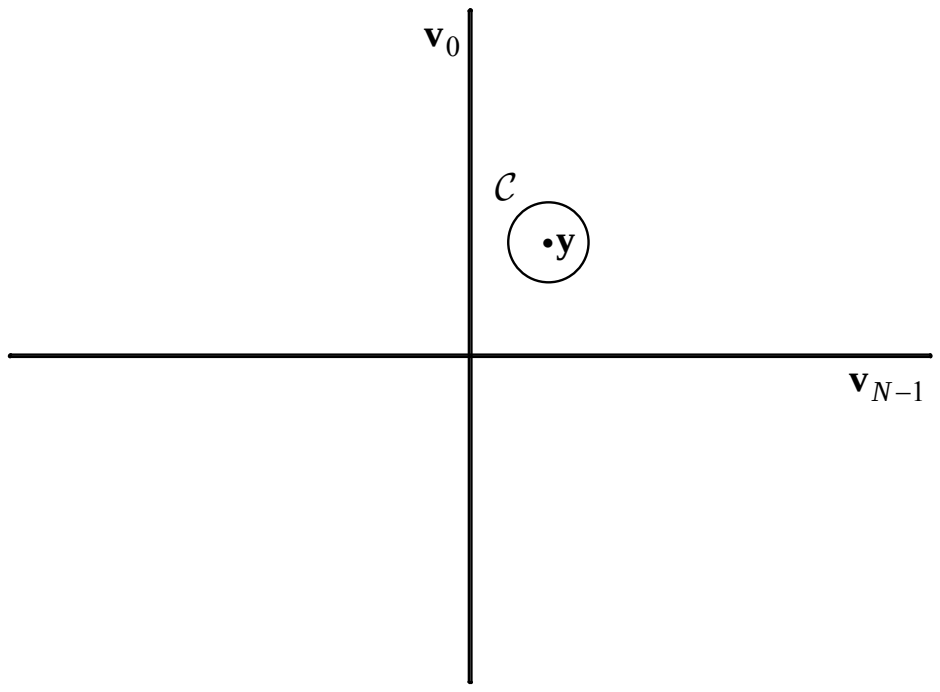

(a)

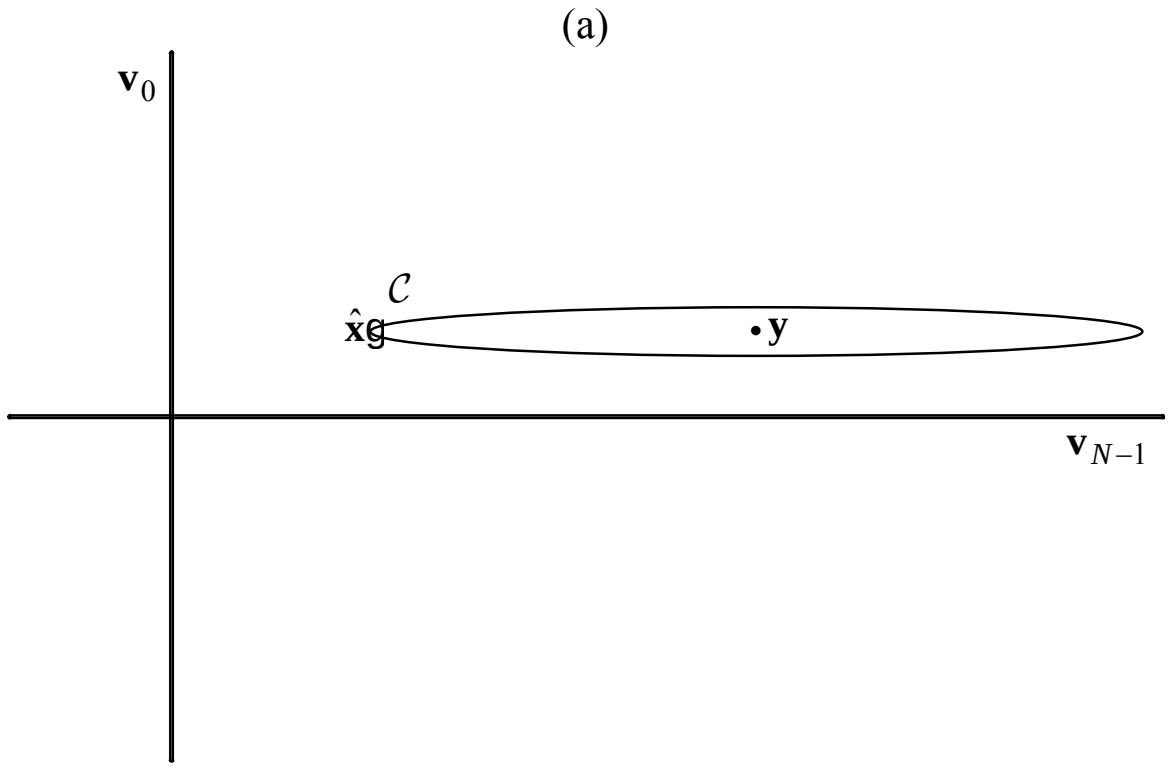

(b)

Fig. 1. Illustration of minimax robust reconstruction geometry – (a) with respect to $\|\cdot\|_2$; (b) with respect to $\|\cdot\|_{\mathbf{K}}$.

It is now easy to see that the stability of the solution relative to $\mathbf{y}$ depends greatly on the location of $\mathbf{y}$ and the shape of $\mathcal{C}$. For example, in the situation illustrated in Figure 1(b), a small change in the coordinate of the nominal observation vector in the positive $\mathbf{v}_{N-1}$ direction will produce a relatively large change in $\hat{\mathbf{x}}$, particularly if $\mathbf{K}$ is ill-conditioned. This is illustrated in Figure 2(a), which indicates a large change in $\hat{\mathbf{x}}$ from the value illustrated in Figure 1(b). On the other hand, a small change of the coordinate in any other direction, for example, even in the negative $\mathbf{v}_{N-1}$ direction, will have much less impact on $\hat{\mathbf{x}}$. This is illustrated in Figure 2(b), which now indicates a much smaller change in $\hat{\mathbf{x}}$ from the value illustrated in Figure 1(b). Hence, even if $\mathbf{K}$ is ill-conditioned, the minimax robust reconstruction will be much less sensitive to small changes in observation *on the average* than interpolation.

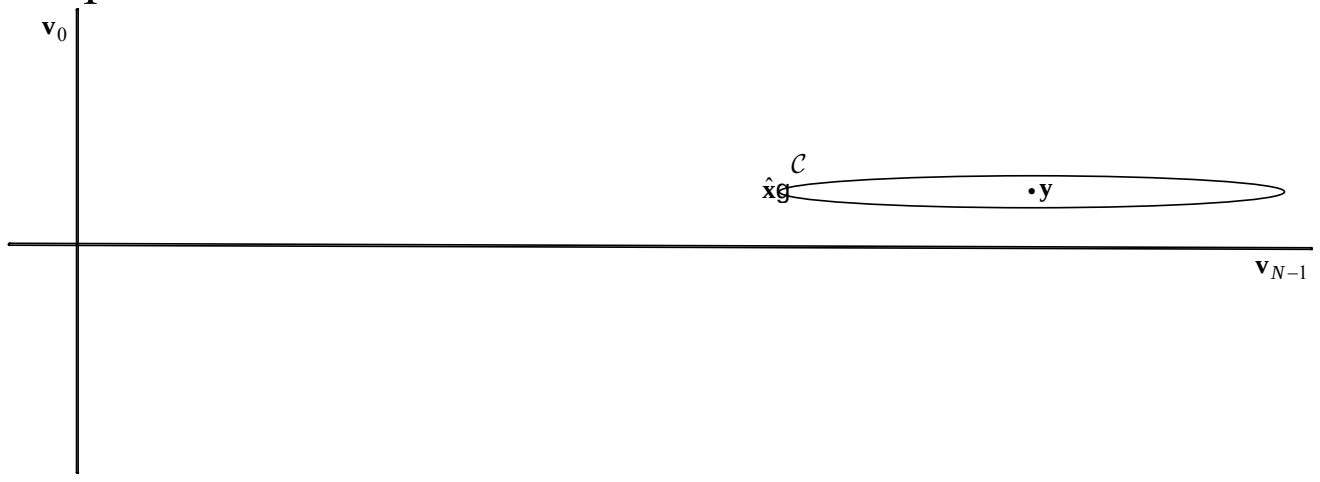

(a)

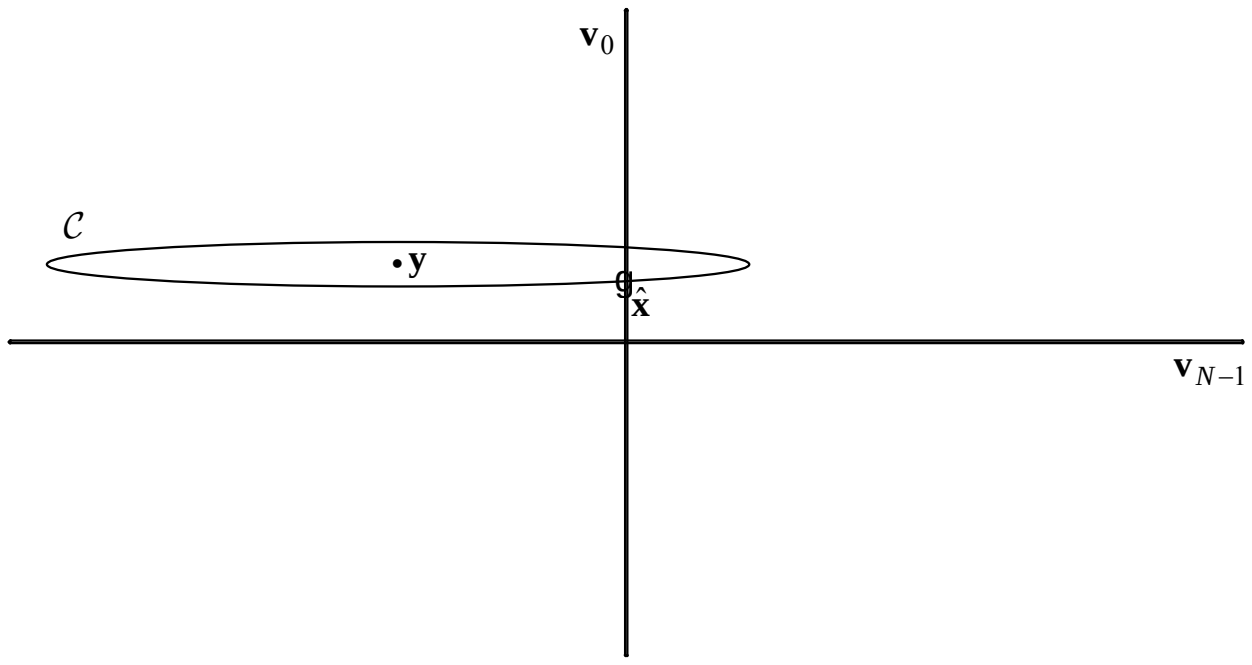

(b)

Fig. 2. Illustration of stability properties of minimax robust reconstruction – (a) small change in the positive $\mathbf{v}_{N-1}$ direction; (b) small change in the negative $\mathbf{v}_{N-1}$ direction.

In fact, minimax robust reconstruction in general can be made unconditionally stable if the set $\mathcal{C}$ is chosen appropriately. For example, suppose $\mathcal{C} = \mathcal{C}_s \times \mathcal{C}_u$ ,where $\mathcal{C}_s$ is a convex subset of the subspace spanned by the "stable" eigenvectors of **K,** and $\mathcal{C}_u$ is a convex subset of the subspace spanned by the "unstable" eigenvectors of **K.**[5] If the origin is contained in the interior of $\mathcal{C}_u$ but is not contained in $\mathcal{C}_s$, then Problem (1) will have a non-trivial ($\hat{\mathbf{x}} \neq 0$) but unconditionally stable solution. A simple example of this is illustrated in Figure 3.

For this figure, we have again illustrated the two-dimensional subspace spanned by $\mathbf{v}_0$ and $\mathbf{v}_{N-1}$, but in this case we let $\mathcal{C} = \mathcal{C}_s \times \mathcal{C}_u$, where $\mathcal{C}_s$ and $\mathcal{C}_u$ are arbitrary convex sets, which in one dimension just correspond to real-valued intervals. Hence, for this example, the uncertainty region becomes a simple rectangle that is much more elongated with respect to the $\|\cdot\|_{\mathbf{K}}$ norm than with respect to the $\|\cdot\|_2$ norm.

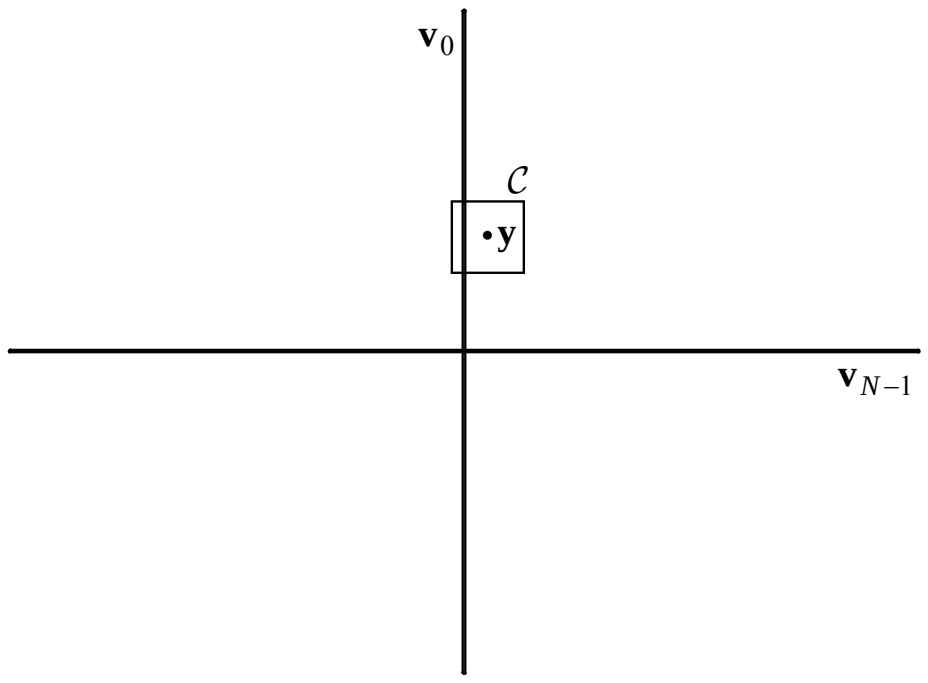

(a)

[5] Here, the terms stable and unstable refer to the relative magnitudes of the eigenvalues. That is, the stable eigenvalues are those with relative magnitude larger than some threshold and the unstable eigenvalues are those with relative magnitude smaller than the same threshold. This is analogous to the manner in which the singular values would be chosen in the singular value decomposition of **K.**

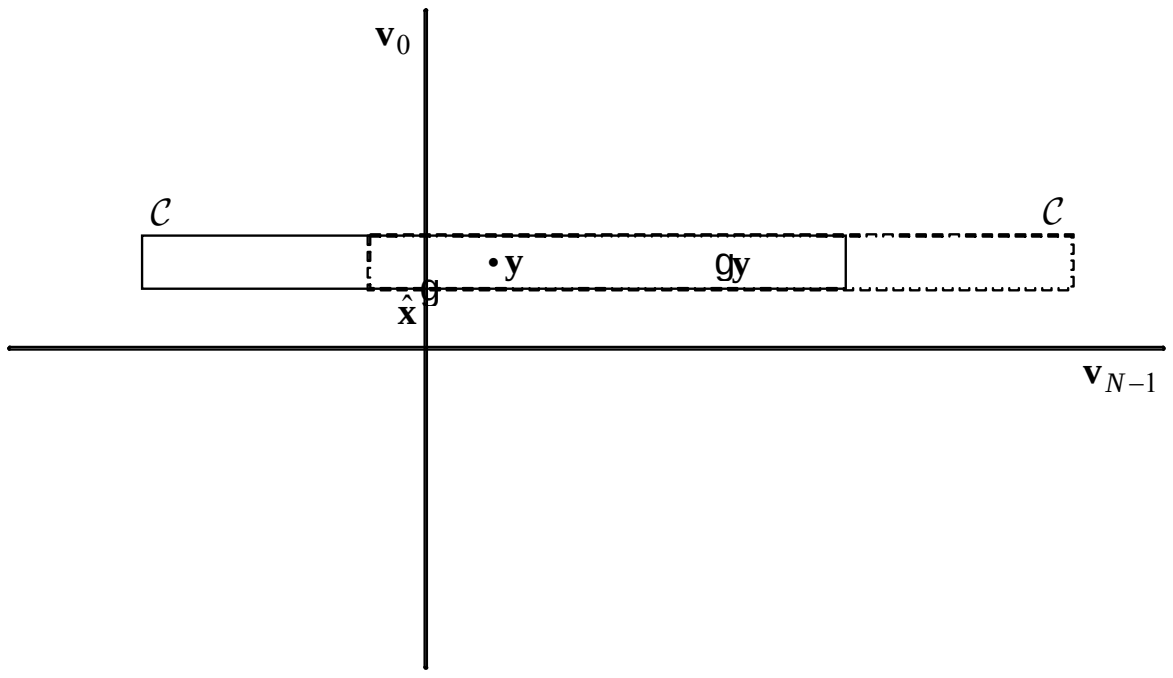

(b)

Fig. 3. Illustration of unconditionally stable minimax robust reconstruction geometry – (a) with respect to $\|\cdot\|_2$ ; (b) with respect to $\|\cdot\|_{\mathbf{K}}$ (solid line indicates original location of $\mathcal{C}$, dashed line indicates position after small change in coordinate in the $\mathbf{v}_{N-1}$ direction).

As Figure 3 indicates, since the origin is contained in the interval representing $\mathcal{C}_u$ , a small change in the coordinates of the nominal observed vector $\mathbf{y}$ in the $\mathbf{v}_{N-1}$ direction will produce no change in the coordinate of the minimax reconstruction $\hat{\mathbf{x}}$ in that direction, which is zero in both cases.

Interestingly, constructing $\mathcal{C}$ in this fashion is analogous to using the Moore-Penrose pseudo-inverse [26] of $\mathbf{K}$ rather than the true inverse of $\mathbf{K}$ when finding a solution in the case of minimum-norm interpolation. That is, in both cases, the projection of the solution onto the unstable eigenvectors is zero. In fact, if the set $\mathbb{R}^{\mathcal{O}}$ contains only the point representing the projection of the nominal observed vector $\mathbf{y}$ onto the span of the stable eigenvectors, then this approach is identical to using the pseudo-inverse of $\mathbf{K}$ to solve the minimum-norm interpolation problem. It follows that minimax robust reconstruction can also be regarded as a generalization of regularization applied to minimum-norm interpolation. Here again, the equivalence established between minimax robust reconstruction and matrix regularization reveals the previously unrecognized minimax optimality property of a function approximation technique that is generally developed primarily in a heuristic context.

## IV. Conclusion

In this paper, we have introduced and studied a multidimensional RKHS function approximation technique that can be viewed as a generalization of scattered data interpolation. The desired approximation is derived by solving a minimax problem with respect to a uniformly bounded uncertainty set of admissible functions that are required to satisfy a known set of convex constraints on the observations. We refer to this approach as minimax robust reconstruction.

We have demonstrated that this approach to function approximation has several significant advantages. From a practical point of view, a major advantage is that either empirical or *a-priori* information about the statistical distribution of noise in the observations can be easily and explicitly incorporated in the selection of the uncertainty set $\mathcal{U}$. In addition, this approach places the function approximation problem in a very geometrical context that facilitates the choice of the uncertainty set and provides insight into sensitivity and error analysis. From a theoretical perspective, minimax robust reconstruction can be regarded as a unified approach to function approximation that includes many other popular techniques as special cases and establishes the previously unrecognized minimax optimality associated with many of these techniques.

The work presented here can be extended easily to handle both approximation in a general Hilbert space and vector-valued functions in a very similar fashion. In addition, the approach can be extended to identify minimax robust reconstructions in situations where the appropriate norm on the function space is itself not completely identified. This might be useful, for example, in situations where the observed functions are regarded as realizations of a Gaussian random process for which the covariance function is known only to belong to a convex set of possible covariance functions. Extensions such as these, as well as other questions of interest regarding the proposed approach, will be studied in future work.

**Richard J. Barton** received a B.A. degree in actuarial science and finance in 1976, an M.S. in Mathematics in 1984, and a Ph.D. in Electrical Engineering in 1989, all from the University of Illinois in Urbana-Champaign.

From 1989 to 1997, he worked for ORINCON Corporation in San Diego, CA. From 1997 to 2006, he was in academia, first in the Electrical and Computer Engineering Department at Iowa State University in Ames, IA and then in the Electrical and Computer Engineering Department at University of Houston in Houston, TX. In 2006, he moved to NASA Johnson Space Center in Houston, where he is currently employed in the Electromagnetic Systems Branch.

Dr. Barton's research interests span many different aspects of statistical signal processing. Past research contributions have been in the areas of robust signal detection and estimation, signal detection in the presence of long-term dependent noise, and applications of wavelet transforms and higher-order statistics to pattern recognition and signal classification. His current research interests include location estimation in wireless environments, detection and estimation using surface acoustic wave (SAW) RFID tags, cooperative communication techniques, and signal processing in low-power wireless sensor networks.